\documentclass[11pt]{article}
\usepackage{amsmath}
\usepackage{pxfonts}
\usepackage{yfonts}
\usepackage{dsfont}
\usepackage{marvosym}
\usepackage{graphicx}
\usepackage{relsize}

\def\bel{\begin{equation}\label}
\def\eeq{\end{equation}}
\def\ds{\displaystyle}
\def\endproof{\hphantom{MM}
\hfill\llap{$\square$}\goodbreak}

\def\mt{\longrightarrow}
\def\v{\vskip 1em}
\def\vsk{\vskip 40em}

\def\R{\mathds R}
\def\Z{\mathds Z}
\def\C{\mathfrak{C}}
\def\Cx{\mathds C}

\def\Re{{\bf Re}}
\def\Im {{\bf Im}}

\def\F{\mathfrak{F}}

\def\J{{\bf J}}

\def\L{{\bf L}}

\def\p{{\partial}}

\def\i{{\bf i}}

\def\Hat{\widehat}

\def\I{{\bf I}}

\def\M{{\bf M}}

\def\Rec{{\bf R}}

\def\alpha{\alphaup}
\def\beta{\betaup}
\def\gamma{\gammaup}
\def\delta{\deltaup}
\def\xi{{\xiup}}
\def\eta{{\etaup}}
\def\tau{{\tauup}}
\def\rho{{\rhoup}}
\def\phi{{\phiup}}
\def\psi{{\psiup}}
\def\lambda{{\lambdaup}}
\def\omega{\omegaup}
\def\varphi{{\varphiup}}
\def\gamma{{\gammaup}}
\def\c{\mathfrak{c}}

\def\e{{\bf e}}
\def\u{{\bf u}}
\def\p{{\bf v}}

\begin{document}
 \[\begin{array}{cc}\hbox{\LARGE{\bf Fractional integration with singularity on unit sphere }}\end{array}\] 
 
  \[\hbox{Zipeng Wang}\]
 \begin{abstract}
We study a family of convolution operators whose kernels have  a singularity on the unit sphere $\mathds{S}^{n-1}$. As a result, we prove the regarding $\L^p\mt\L^q$-Sobolev inequalities.
\end{abstract}
 \section{Introduction}
 \setcounter{equation}{0}
In this paper, we revisit on a classical problem  that has been previously investigated by Strichartz \cite{Strichartz}. 

Let 
$\lambda(\alpha)\doteq{n+1\over 2}\left(1-{\alpha\over n}\right),\alpha\in\Cx$.
For $\Re\lambda(\alpha)<1$, define 
\bel{Omega^alpha}
\begin{array}{ccc}\ds
\Omega^{\alpha}(x)~\doteq~\left\{\begin{array}{lr}\ds
\pi^{-\lambda(\alpha)}\Gamma^{-1}\left(1-\lambda(\alpha)\right)  \left({1\over 1-|x|^2}\right)^{\lambda(\alpha)}, \qquad |x|<1,
\\\\ \ds~~~~~~~~~~~~~~~~~~~~~~
0,\qquad\qquad~~~~~~~~~~~~~~~~~~~~~~~~|x|\ge1
\end{array}\right.
\end{array}
\eeq
whose Fourier transform equals
 \bel{Omega^alpha Transform} 
\begin{array}{lr}\ds
\Hat{\Omega}^{\alpha}(\xi)~=~\left({1\over|\xi|}\right)^{{n\over 2}-\lambda(\alpha)} \J_{{n\over 2}-\lambda(\alpha)}\Big(2\pi|\xi|\Big).
\end{array}
\eeq
$\Gamma$ and $\J$ denote for Gamma and Bessel functions. See chapter IV in the book by Stein and Weiss \cite{Stein-Weiss*}. 

Observe that $\Omega^\alpha$ in (\ref{Omega^alpha}) has an non-integrable singularity if $\Re\lambda(\alpha)\ge1$. However, we can extensively define  $\Omega^\alpha$ for every $\alpha\in\Cx$ by its Fourier transform in (\ref{Omega^alpha Transform}). 

$\diamond$ {\small $\L^p_s$ is the generalized Sobolev space for $s>0$ and $1<p<\infty$}.

$\diamond$ {\small Throughout, we regard $\C$ as a generic constant depending on its sub-indices}.

{\bf Theorem One}~~{\it Let $\Omega^\alpha$ defined by its Fourier transform in (\ref{Omega^alpha Transform}) for $0<\alpha<n$.  We have
\bel{RESULT ONE}
\left\| f\ast\Omega^{\alpha}\right\|_{\L^q(\R^n)}~\leq~\C_{p~q}~\left\| f\right\|_{\L^p_s(\R^n)},\qquad s>0
\eeq
\bel{FORMULA ONE}
\begin{array}{cc}\ds
\hbox{if}\qquad {\alpha\over n}~=~{1\over p}-{1\over q}\qquad \hbox{and}
\\\\ \ds
{n-1\over 2n-2+4s}+\left({4s+n-1\over 2n-2+4s}\right){\alpha\over n}~<~{1\over p}~<~{n-1+4s\over 2n-2+4s}+\left({n-1\over 2n-2+4s}\right){\alpha\over n}.
\end{array}
\eeq}

In particular, at $s={1\over 2}$, we have
\bel{Result 1/2 p q}
\begin{array}{cc}\ds
\left\| f\ast \Omega^\alpha\right\|_{\L^q(\R^n)}~\leq~\C_{p~q}~\left\| f\right\|_{\L^p_{1/ 2}(\R^n)}
\\\\ \ds
\hbox{if}\qquad 
{\alpha\over n}~=~{1\over p}-{1\over q},\qquad {n-1\over 2n}+\left({n+1\over 2n}\right){\alpha\over n}~<~{1\over p}~<~{n+1\over 2n}+\left({n-1\over 2n}\right){\alpha\over n}.
\end{array}
\eeq
 \begin{figure}[h]
\centering
\includegraphics[scale=0.34]{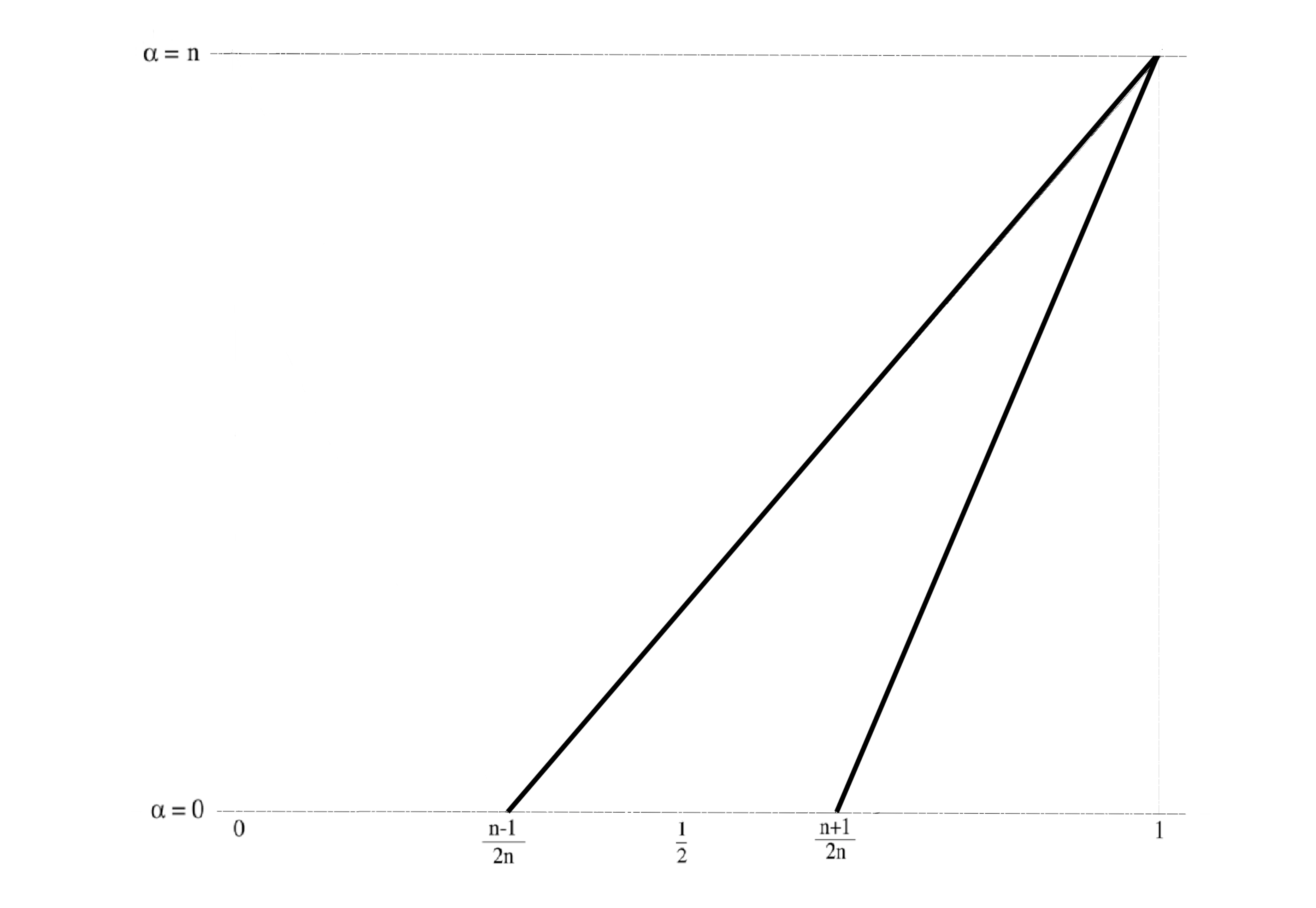}
\end{figure} 

One application of {\bf Theorem One} gives  an {\it a priori} estimate of the wave equation
\bel{Wave}
\begin{array}{cc}\ds
\partial^2_t \u (x,t)-\Delta_x \u(x,t)~=~f(x,t),\qquad (x,t)\in\R^n\times[0,\infty),
\\\\ \ds
\u(x,0)~=~\partial_t\u(x,0)~=~0.
\end{array}
\eeq
{\bf Corollary One}~~ {\it Let $\u$ be the solution of (\ref{Wave}). We have
\bel{u norm}
\begin{array}{ccc}\ds
\left\| \u(\cdot,t)\right\|_{\L^q(\R^n)}~\leq~\C_{p~q}~\int_0^t \left\| f(\cdot, t-r)\right\|_{\L^p_s(\R^n)} r^{-{n-1\over n+1}} dr,\qquad s>0,\qquad t\in[0,\infty)
\\\\ \ds
\hbox{if}\qquad {2\over n+1}~=~{1\over p}-{1\over q}\qquad \hbox{and}
\\\\ \ds
{n-1\over 2n-2+4s}+\left({4s+n-1\over 2n-2+4s}\right){2\over n+1}~<~{1\over p}~<~{n-1+4s\over 2n-2+4s}+\left({n-1\over 2n-2+4s}\right){2\over n+1}.
\end{array}
\eeq}

For more discussion of Sobolev inequalities regarding to wave equations, we refer to the references by Sogge \cite{Sogge}-\cite{Sogge'}. 

{\bf Proof}~~As shown by Strichartz \cite{Strichartz'}, it is well known that
\bel{u Transform}
\Hat{\u}(\xi,t)~=~\int_0^t {\sin[(t-s)|\xi|]\over |\xi|} \Hat{f}(\xi,s)ds.
\eeq
From (\ref{Omega^alpha Transform}), we find $\Omega^{2n\over n+1}(\xi)={1\over \pi}{\sin 2\pi|\xi|\over|\xi|}$. Let $f\in\L^p_s(\R^n)$. By applying {\bf Theorem One}, we have
\bel{u norm Est}
\begin{array}{lr}\ds
\left\| \u(\cdot,t)\right\|_{\L^q(\R^n)}~\leq~\C~\int_0^t \left\| f(\cdot, t-r)\ast\Omega^{2n\over n+1}_r\right\|_{\L^q(\R^n)} rdr
\\\\ \ds~~~~~~~~~~~~~~~~~~~~~
~\leq~\C_{p~q}~\int_0^t \left\| f(\cdot, t-r)\right\|_{\L^p(\R^n)} r^{1-{2n\over n+1}}dr.
\end{array}
\eeq 
\endproof

{\bf Theorem Two}~~{\it Let $\Omega^\alpha$ defined by its Fourier transform in (\ref{Omega^alpha Transform}) for $0<\alpha<\left({n-1\over n+1}\right)n$.  We have
\bel{RESULT TWO}
\begin{array}{ccc}\ds
\left\| f\ast\Omega^\alpha\right\|_{\L^p(\R^n)}
 ~\leq~\C_{p~\alpha}~\left\| f\right\|_{\L^p_s(\R^n)},\qquad s>0
 \\\\ \ds
 \hbox{for}\qquad {n-1\over 2n-2+4s}-\left({n+1\over 2n-2+4s}\right){\alpha\over n}~<~{1\over p}~<~{n-1+4s\over 2n-2+4s}+\left({n+1\over 2n-2+4s}\right){\alpha\over n}.
 \end{array}
\eeq
}

{\bf Remark One}~~{\it For $\left({n-1\over n+1}\right)n\leq\alpha<n$, we have
\bel{L^p-result}
\left\| f\ast \Omega^\alpha\right\|_{\L^p(\R^n)}~\leq~\C_{p~\alpha}~\left\| f\right\|_{\L^p(\R^n)},\qquad 1<p<\infty.
\eeq}
For $s={1\over 2}$ in (\ref{RESULT TWO}), we have
\bel{Result 1/2 p p}
\begin{array}{cc}\ds
\left\| f\ast \Omega^\alpha\right\|_{\L^p(\R^n)}~\leq~\C_{p~\alpha}~\left\| f\right\|_{\L^p_{1/ 2}(\R^n)}
\\\\ \ds
\hbox{for}\qquad 
 {n-1\over 2n}-\left({n+1\over 2n}\right){\alpha\over n}~<~{1\over p}~<~{n+1\over 2n}+\left({n+1\over 2n}\right){\alpha\over n}.
\end{array}
\eeq

\section{Principal Lemmata}
\setcounter{equation}{0}
In order to obtain the desired regularity of $f\ast\Omega^\alpha$, a key estimate is to study  \bel{Omega^alpha natural}
\begin{array}{ccc}\ds
 {^\natural}\Omega^{\alpha}(x)~\doteq~\left\{\begin{array}{lr}\ds
 \left({1\over 1-|x|^2}\right)^{1-{\alpha\over n}}, \qquad |x|<1,
\\\\ \ds~~~~~~~
0,\qquad~~~~~~~~~~~~~|x|\ge1.
\end{array}\right.
\end{array}
\eeq
 
{\bf Lemma One}~~ {\it Let ${^\natural}\Omega^\alpha$ defined by (\ref{Omega^alpha natural}) for $0<\alpha< n$. We have
\bel{Result One}
\begin{array}{cc}\ds
\left\| f\ast {^\natural}\Omega^\alpha\right\|_{\L^q(\R^n)}~\leq~\C_{\alpha~p~q}~\left\| f\right\|_{\L^p(\R^n)},\qquad 1<p\leq q<\infty
\\\\ \ds
\hbox{if}\qquad {\alpha\over n}~\ge~{1\over p}-{1\over q}.
\end{array}
\eeq}

Let $\delta(\alpha)\doteq1-\left({n+1\over 2n}\right)\alpha$ for $0<\Re\alpha<{2n\over n+1}$. 
Define 
\bel{Omega flat}
{^\flat}\Omega^{\alpha}(x)~\doteq~\left\{\begin{array}{lr}\ds
\hbox{\small{$\pi^{-\delta(\alpha)}  \Gamma^{-1}\left(1-\delta(\alpha)\right)$}}
\left({1\over 1-|x|^2}\right)^{\delta(\alpha)},\qquad~~0<|x|<1,
\\\\ \ds~~~~~~~~~~~~~~~~~~~~~~~~
0,\qquad\qquad\qquad\qquad~~~~~~~|x|\ge1
\end{array}
\right.
\eeq
whose Fourier transform equals
\bel{Omega flat Transform}
\begin{array}{lr}\ds
{^\flat}\Hat{\Omega}^{\alpha}(\xi)~=~  \left({1\over|\xi|}\right)^{{n\over 2}-\delta(\alpha)} \J_{{n\over 2}-\delta(\alpha)}\Big(2\pi|\xi|\Big)
\\\\ \ds~~~~~~~~~~~
~=~\left({1\over|\xi|}\right)^{\left({n+1\over 2n}\right)\alpha-{1\over 2}+{n-1\over 2}} 
\J_{\left({n+1\over 2n}\right)\alpha-{1\over 2}+{n-1\over 2}}\Big(2\pi|\xi|\Big).
\end{array}
\eeq
$\diamond$ {\small In the remaining paragraph, we write $\hbox{\bf c}>0$ for some large constant.}

Observe that 
\bel{Kernel compara}
\left|{^\flat}\Omega^{\alpha}(x)\right|~\leq~\C_\alpha ~e^{\hbox{\small{\bf c}}|\Im\alpha|}~{^\natural}\Omega^{\left({n+1\over 2}\right)\Re\alpha}(x)
\eeq
where ${^\natural}\Omega^\alpha$ is defined in (\ref{Omega^alpha natural}).

By  using (\ref{Kernel compara}) and applying {\bf Lemma One}, we have
\bel{Omega flat regularity}
\begin{array}{cc}\ds
\left\| f\ast{^\flat}\Omega^{\alpha}\right\|_{\L^q(\R^n)}~\leq~\C_{\alpha~p~q}~e^{\hbox{\bf \small{c}}|\Im\alpha|}~\left\|f\right\|_{\L^p(\R^n)},\qquad 1<p\leq q<\infty
\\\\ \ds
\hbox{if}\qquad \left({n+1\over 2n}\right)\Re\alpha~\ge~{1\over p}-{1\over q}.
\end{array}
\eeq
For $s\ge0$ and $\Re\lambda(\alpha)+s<1$, define
\bel{Omega s}
\begin{array}{ccc}\ds
{^s}\Omega^{\alpha}(x)~\doteq~\left\{\begin{array}{lr}\ds
\pi^{-\lambda(\alpha)-s}\Gamma^{-1}\left(1-\lambda(\alpha)-s\right)  \left({1\over 1-|x|^2}\right)^{\lambda(\alpha)+s}, \qquad |x|<1,
\\\\ \ds~~~~~~~~~~~~~~~~~~~~~~~~~~
0,\qquad\qquad~~~~~~~~~~~~~~~~~~~~~~~~~~~~~~~~|x|\ge1
\end{array}\right.
\end{array}
\eeq
whose Fourier transform equals
\bel{Omega s Transform}
\begin{array}{lr}\ds
{^s}\Hat{\Omega}^{\alpha}(\xi)~=~\left({1\over|\xi|}\right)^{{n\over 2}-\lambda(\alpha)-s} \J_{{n\over 2}-\lambda(\alpha)-s}\Big(2\pi|\xi|\Big)
\\\\ \ds~~~~~~~~~~
~=~\left({1\over|\xi|}\right)^{\left({n+1\over 2n}\right)\alpha-{1\over 2}-s} \J_{\left({n+1\over 2n}\right)\alpha-{1\over 2}-s}\Big(2\pi|\xi|\Big).
\end{array}
\eeq
Let $\omega_s$ defined by $\Hat{\omega}_s(\xi)\doteq\left({1\over 1+|\xi|^2}\right)^{s\over 2}$ for $s\ge0$. Note that $\Hat{\omega}_s, s\ge0$ is a $\L^p$-Fourier multiplier. Namely, we have
\bel{L^p omega}
\left\| f\ast\omega_s\right\|_{\L^p(\R^n)}~\leq~\C_p~\left\| f\right\|_{\L^p(\R^n)},\qquad 1<p<\infty.
\eeq
Regarding estimates can be found in chapter VI of  Stein \cite{Stein}. 
\vsk

{\bf Lemma Two}~~{\it Let ${^s}\Omega^\alpha$ defined by its Fourier transform in (\ref{Omega s Transform}) for $0<\Re\alpha<n$. 

For $0\leq s\leq{1\over 2}$, we have
 \bel{Result Two < EST}
  \begin{array}{cc}\ds
 \left\| f\ast\omega_s\ast{^s}\Omega^{\alpha}\right\|_{\L^q(\R^n)}~\leq~\C_{p~q~s}~e^{\hbox{\small{\bf c}}|\Im\alpha|}~\left\| f\right\|_{\L^p(\R^n)}
 \\\\ \ds
 \hbox{ if}\qquad
 (1-s){\Re\alpha\over n}~=~{1\over p}-{1\over q},\qquad {1\over 2}+{\Re\alpha\over n}\left({1\over 2}-s\right)~\leq~{1\over p}~\leq~{1\over 2}+{\Re\alpha\over 2n}.
 \end{array}
 \eeq
 For $s\ge{1\over 2}$, we have
 \bel{Result Two > EST}
  \begin{array}{cc}\ds
 \left\| f\ast\omega_s\ast{^s}\Omega^{\alpha}\right\|_{\L^q(\R^n)}~\leq~\C_{p~q~s}~e^{\hbox{\small{\bf c}}|\Im\alpha|}~\left\| f\right\|_{\L^p(\R^n)}
 \\\\ \ds
 \hbox{ if}\qquad
{\Re\alpha\over 2n}~=~{1\over p}-{1\over q},\qquad {1\over 2}~\leq~{1\over p}~\leq~{1\over 2}+{\Re\alpha\over 2n}.
 \end{array}
 \eeq 
For every $s\ge0$, we have
\bel{Result Two L^2}
\left\| f\ast\omega_s\ast{^s}\Omega^{\alpha}\right\|_{\L^2(\R^n)}~\leq~\C_{s}~e^{\hbox{\small{\bf c}}|\Im\alpha|}~\left\| f\right\|_{\L^2(\R^n)}.
\eeq}

We prove {\bf Lemma One} in section 3 and {\bf Lemma Two} in section 4. We finish the proof of {\bf Theorem One} and {\bf Theorem Two} in section 5.

\section{Proof of Lemma One}
\setcounter{equation}{0}

We prove the lemma in the same sprit of Hedberg \cite{Hedberg}. 
Let $\I^\natural_\alpha f\doteq f\ast{^\natural}\Omega^\alpha$ and assume $f\ge0$. For every  $0\leq\ell\in\Z$,  define
 \bel{Partial}
 \begin{array}{cc}\ds
\Big( \Delta_\ell \I_\alpha^\natural f\Big)(x)~\doteq~\int_{\mathcal{S}_\ell} f(x-u)\left({1\over 1-|u|^2 }\right)^{1-{\alpha\over n}}du,
\\\\ \ds
\mathcal{S}_\ell~\doteq~\left\{u\in\R^n~\colon 2^{-\ell-1}\leq 1-|u| < 2^{-\ell}\right\}.
\end{array}
\eeq 
Momentarily, we consider the partial sum operator 
\bel{Partial sum}
\I_\alpha^\rho~\doteq~\sum_{0\leq\ell\leq \rho} \Delta_\ell\I_\alpha^\natural\qquad\hbox{for $\rho$ sufficiently large.}
\eeq
Let $\left\{u^\p_\rho\right\}_\p$ be
a collection of points that are  equally distributed on  $\mathds{S}^{n-1}\subset\R^n$ with grid length equal to $2^{-\rho}$ multiplied by a suitable constant.   
Define the {\it narrow cone} 
 \bel{cone}
\Gamma^\p_\rho~\doteq~\left\{u\in\R^n~\colon~\left| {u\over |u|}-u^\p_\rho\right|~\leq~2^{-\rho}\right\}.
\eeq
Let $\Rec^\p_{\rho~\ell}( u^\p_\rho)$ be the rectangle  centered on $u^\p_\rho$ with one side parallel to $u^\p_\rho$ of side length  $5\times2^{-\ell}$    and others  perpendicular to $u^\p_\rho$ equal to $5\times2^{-\rho}$. 

Observe that
\bel{Inclusion R^v_j} 
\begin{array}{cc}\ds
\mathcal{S}_\ell\cap \Gamma^\p_\rho~~\subset~~\Rec_{\rho~\ell}^\p ( u^\p_\rho),\qquad 0\leq\ell\leq \rho.
\end{array}
\eeq
 \begin{figure}[h]
\centering
\includegraphics[scale=0.40]{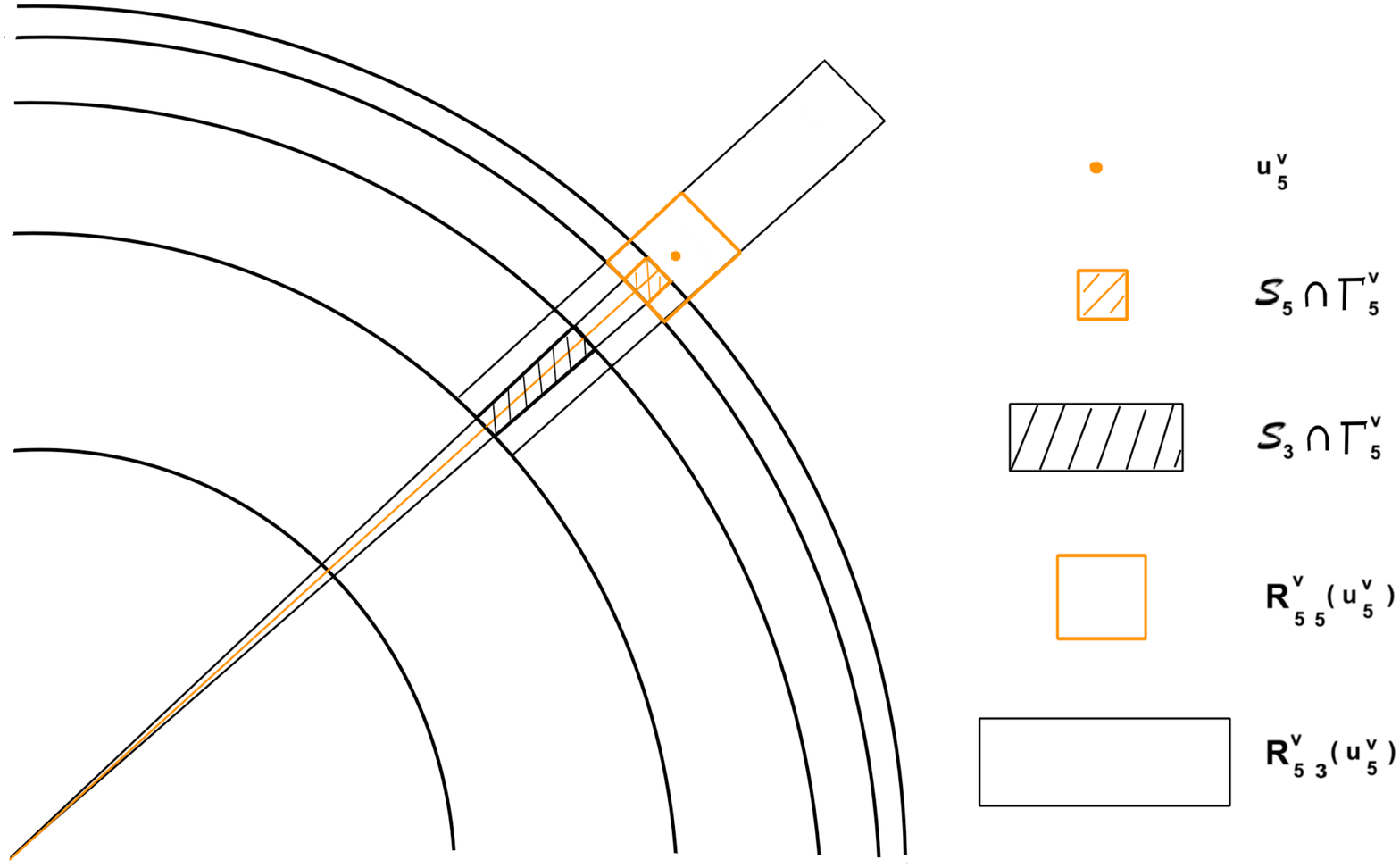}
\caption{\small{$\rho=5$ and $\ell=3,5$.}}
\end{figure}

Let
\bel{R^v_j o}
\Rec^\p_{\rho~\ell}~\doteq~\left\{u\in\R^n\colon u+u^\p_\rho\in \Rec^\p_{\rho~\ell }( u^\p_\rho)\right\}
\eeq
where 
\bel{R^v_j size}
\left| \Rec^\p_{\rho~\ell }( u^\p_\rho)\right|~=~5^n ~ 2^{-\ell}2^{-\rho(n-1)}.
\eeq
We define
\bel{M^v} 
\begin{array}{lr}\ds
\Big( \M^\p_\rho f\Big)(x)~\doteq~\sup_{0\leq\ell\leq \rho} ~{1\over 2^{-\ell}2^{-\rho(n-1)}}\int_{\mathcal{S}_\ell\cap\Gamma^\p_\rho} f(x-u)du
 \\\\ \ds~~~~~~~~~~~~~~~
 ~\leq~\sup_{0\leq\ell\leq \rho} ~{5^n\over |\Rec^\p_{\rho~\ell}|}\int_{\Rec_{\rho~\ell}^\p} f\left(x-u^\p_\rho-u\right)du.
 \end{array}
 \eeq
 Such  a maximal operator  defined $w.r.t$ all rectangles having a same degree is  bounded on $\L^p(\R^n)$
  for $1<p<\infty$.   
  
 There are at most $\C 2^{\rho\left(n-1\right)}$ elements  
 in $\left\{u^\p_\rho\right\}_\p$. Consider
\bel{M*}
\Big(\M_\rho f\Big)(x) ~\doteq~2^{-\rho(n-1)}\sum_\p \Big(\M^\p_\rho f\Big)(x).
\eeq
{\bf Remark Two} {\it $\M_\rho$ in (\ref{M*}) is NOT the strong maximal operator defined $w.r.t$ all rectangles with arbitrary directions. The later has unboundedness on  every $\L^p$-space for $1\leq p<\infty$ due to the Besicovitch construction. See Fefferman \cite{Fefferman} or chapter X and XI in the book by Stein \cite{Stein}.}

By using (\ref{M*}) and applying Minkowski inequality, we have 
\bel{L^p M}
\begin{array}{lr}\ds
\left\|\M_\rho f\right\|_{\L^p(\R^n)}
~\leq~\C_p~\left\|f\right\|_{\L^p(\R^n)},\qquad 1<p<\infty.
\end{array}
\eeq
In particular, (\ref{L^p M}) is independent from $\rho$.

Let $f\in\L^p(\R^n)$. We claim 
  \bel{Regularity est}
 \begin{array}{cc}\ds
 \Big( \I_\alpha^\rho f\Big)(x)
 ~\leq~\C_{p~q}~  \Big(\M_\rho f\Big)^{p\over q}(x) \left\| f\right\|_{\L^p(\R^n)}^{1-{p\over q}}
\end{array}
\eeq
for
\bel{formula}
{\alpha\over n}~=~{1\over p}-{1\over q},  \qquad1<p<q<\infty.
\eeq
By using the $\L^p$-boundedness of $\M_\rho$ in (\ref{L^p M}), we  find 
\bel{Norm Ineq natural eta}
\begin{array}{lr}\ds
\left\|\I_\alpha^\rho f\right\|_{\L^q(\R^n)}~\leq~\C_{p~q}~\left\{\int_{\R^n}\Big(\M_\rho f\Big)^{p}(x)dx\right\}^{1\over q} \left\| f\right\|_{\L^p(\R^n)}^{1-{p\over q}}
\\\\ \ds~~~~~~~~~~~~~~~~~~~
~\leq~\C_{p~q}~\left\|f\right\|_{\L^p(\R^n)}.
\end{array}
\eeq
By allowing $\rho\mt\infty$ in (\ref{Norm Ineq natural eta}) and using Fatou's lemma together with the monotone convergence theorem, we conclude
\bel{Norm Ineq natural}
\begin{array}{lr}\ds
\left\|\I_\alpha^\natural f\right\|_{\L^q(\R^n)}~\leq~\C_{p~q}~\left\|f\right\|_{\L^p(\R^n)}.
\end{array}
\eeq
On the other hand, from (\ref{Omega^alpha natural}), it is clear that $\I_\alpha^\natural$ is bounded on $\L^p(\R^n)$ for every $1\leq p\leq\infty$. The desired result in (\ref{Result One}) can be obtained by interpolation.

In order to prove (\ref{Regularity est}), we define $\sigma(x)\in\R$  implicitly by
\bel{sigma}
2^{\sigma(x)/p}~\doteq~{ \left(\M_\rho f\right)(x)\over\left\| f\right\|_{\L^p(\R^n)}}.
\eeq
Let $\mathcal{S}_\ell$ and $\Gamma^\p_\rho$ be given in (\ref{Partial}) and (\ref{cone}). From (\ref{Inclusion R^v_j})-(\ref{R^v_j size}), we have
\bel{Partial Est1 v}
\begin{array}{lr}\ds
\int_{\mathcal{S}_\ell \cap\Gamma^\p_\rho}f(x-u)\left({1\over 1-|u|^2}\right)^{1-{\alpha\over n}} du
~=~\int_{\mathcal{S}_\ell \cap\Gamma^\p_\rho}f(x-u)\left({1\over 1+|u|}\right)^{1-{\alpha\over n}} \left({1\over 1-|u|}\right)^{1-{\alpha\over n}}du
\\\\ \ds~~~~~~~~~~~~~~~~~
~\leq~\C~   2^{\ell(1-{\alpha\over n})}\int_{\mathcal{S}_\ell \cap\Gamma^\p_\rho}f(x-u)du
\\\\ \ds~~~~~~~~~~~~~~~~~
~=~\C~2^{\ell\left(1-{\alpha\over n}\right) }2^{-\ell}2^{-\rho(n-1)}\left\{{1\over 2^{-\ell}2^{-\jmath(n-1)}}\int_{\mathcal{S}_\ell \cap\Gamma^\p_\rho}f(x-u)du\right\}
\\\\ \ds~~~~~~~~~~~~~~~~~
~=~\C~2^{-\ell\left({\alpha\over n}\right) }2^{-\rho(n-1)}\left\{{1\over 2^{-\ell}2^{-\rho(n-1)}}\int_{\mathcal{S}_\ell \cap\Gamma^\p_\rho}f(x-u)du\right\}
\\\\ \ds~~~~~~~~~~~~~~~~~
~\leq~\C~2^{-\ell\left({\alpha\over n}\right) }2^{-\rho(n-1)}\Big(\M_\rho^\p f\Big)(x).
\end{array}
\eeq
Therefore,
\bel{Partial Est1}
\begin{array}{lr}\ds
\int_{\mathcal{S}_\ell }f(x-u)\left({1\over 1-|u|^2}\right)^{1-{\alpha\over n}} du
~\leq~\sum_\p \int_{\mathcal{S}_\ell \cap\Gamma^\p_\rho}f(x-u)\left({1\over 1-|u|^2}\right)^{1-{\alpha\over n}} du
\\\\ \ds~~~~~~~~~~~~~~~~~~~~~~~~~~~~~~~~~~~~~~~~~~~~~~~~
~\leq~\C~2^{-\ell\left({\alpha\over n}\right) }2^{-\rho(n-1)}\sum_\p\Big(\M_\rho^\p f\Big)(x)
\\\\ \ds~~~~~~~~~~~~~~~~~~~~~~~~~~~~~~~~~~~~~~~~~~~~~~~~
~=~\C~2^{-\ell\left({\alpha\over n}\right) }\Big(\M_\rho f\Big)(x).
\end{array}
\eeq
On the other hand, by applying H\"{o}lder inequality, we have
\bel{Partial Est2 Holder}
\begin{array}{lr}\ds
\int_{\mathcal{S}_\ell }f(x-u)\left({1\over 1-|u|^2}\right)^{1-{\alpha\over n}} du
~\leq~\left\| f\right\|_{\L^p(\R^n)}\left\{\int_{\mathcal{S}_\ell }\left({1\over 1-|u|^2}\right)^{\left(1-{\alpha\over n}\right)\left({p\over p-1}\right)} du\right\}^{p-1\over p}
\end{array}
\eeq
where
\bel{Partial Est2}
\begin{array}{lr}\ds
\left\{\int_{\mathcal{S}_\ell }\left({1\over 1-|u|^2}\right)^{\left(1-{\alpha\over n}\right)\left({p\over p-1}\right)} du\right\}^{p-1\over p}
~=~\left\{\int_{\mathcal{S}_\ell }\left({1\over 1+|u|}\right)^{\left(1-{\alpha\over n}\right)\left({p\over p-1}\right)} \left({1\over 1-|u|}\right)^{\left(1-{\alpha\over n}\right)\left({p\over p-1}\right)} du\right\}^{p-1\over p}
\\\\ \ds~~~~~~~~~~~~~~~~~~~~~~~~~~~~~~~~~~~~~~~~~~~~~~~~~~~~
~\leq~\C~\left[   {1\over 2^{-\ell\left(1-{\alpha\over n}\right)\left({p\over p-1}\right)}} 2^{-\ell}\right]^{p-1\over p}
~=~\C~2^{\ell\left(1-{\alpha\over n} \right)}2^{-\ell\left({p-1\over p}\right)}
~=~\C~2^{-\ell \left({\alpha\over n}-{1\over p}\right)}.
\end{array}
\eeq

Consider  $\sigma(x)\leq\ell\leq \rho$.
By inserting (\ref{sigma}) into (\ref{Partial Est1}), we find
\bel{Est1 partial}
\begin{array}{lr}\ds
\int_{\mathcal{S}_\ell}f(x-u)\left({1\over 1-|u|^2}\right)^{1-{\alpha\over n}} du
~\leq~\C~2^{-\ell\left({\alpha\over n}\right)}\Big(\M_\jmath f\Big)(x)
\\\\ \ds~~~~~~~
~=~\C~2^{\left(\sigma(x)-\ell\right)\left({\alpha\over n}\right)} 2^{-\sigma(x)\left({1\over p}-{1\over q}\right)} \Big(\M_\rho f\Big)(x)
\\\\ \ds~~~~~~~
~=~\C~2^{\left(\sigma(x)-\ell\right)\left({\alpha\over n}\right)} \left\{{\left\| f\right\|_{\L^p(\R^n)}\over \left(\M_\rho f\right)(x)}\right\}^{1-{p\over q}}  \Big(\M_\rho f\Big)(x)
\\\\ \ds~~~~~~~
~=~\C~2^{\left(\sigma(x)-\ell\right)\left({\alpha\over n}\right)}\Big(\M_\rho f\Big)^{p\over q}(x)\left\| f\right\|_{\L^p(\R^n)}^{1-{p\over q}}.
\end{array}
\eeq
By summing over all $\sigma(x)\leq\ell\leq \rho$, we have
\bel{Est1 Partial}
\begin{array}{lr}\ds
\sum_{\sigma(x)\leq\ell\leq \rho}\int_{\mathcal{S}_\ell}f(x-u)\left({1\over 1-|u|^2}\right)^{1-{\alpha\over n}} du
~\leq~\C\sum_{\sigma(x)\leq\ell\leq \rho} 2^{\left(\sigma(x)-\ell\right)\left({\alpha\over n}\right)}\Big(\M_\rho f\Big)^{p\over q}(x)\left\| f\right\|_{\L^p(\R^n)}^{1-{p\over q}}
\\\\ \ds~~~~~~~~~~~~~~~~~~~~~~~~~~~~~~~~~~~~~~~~~~~~~~~~~~~~~~~~~~~
~\leq~\C_{p~q}~\Big(\M_\rho f\Big)^{p\over q}(x)\left\| f\right\|_{\L^p(\R^n)}^{1-{p\over q}}.
\end{array}
\eeq
Consider $0\leq\ell<\sigma(x)$.
By inserting (\ref{sigma}) into (\ref{Partial Est2 Holder})-(\ref{Partial Est2}),  we find
\bel{Est2 partial}
\begin{array}{lr}\ds
\int_{\mathcal{S}_\ell }f(x-u)\left({1\over 1-|u|^2}\right)^{1-{\alpha\over n}} du
~\leq~\C~\left\| f\right\|_{\L^p(\R^n)} ~ 2^{-\ell\left({\alpha\over n}-{1\over p}\right)}
\\\\ \ds~~~~~~~
~=~\C~2^{\left(\sigma(x)-\ell\right)\left({\alpha\over n}-{1\over p}\right)} 2^{-\sigma(x)\left({\alpha\over n}-{1\over p}\right)} \left\| f\right\|_{\L^p(\R^n)} 
\\\\ \ds~~~~~~~
~=~\C~2^{-\left(\sigma(x)-\ell\right)\left({1\over q}\right)} 
\left\| f\right\|_{\L^p(\R^n)}\left\{{\left\| f\right\|_{\L^p(\R^n)}\over \left(\M_\rho f\right)(x)}\right\}^{-{p\over q}} 
\\\\ \ds~~~~~~~
~=~\C~2^{-\left(\sigma(x)-\ell\right)\left({1\over q}\right)} \Big(\M_\rho f\Big)^{p\over q}(x) \left\| f\right\|_{\L^p(\R^n)}^{1-{p\over q}}.
\end{array}
\eeq
By summing over all $0\leq\ell<\sigma(x)$, we have
\bel{Est2 Partial}
\begin{array}{lr}\ds
\sum_{0\leq\ell<\sigma(x)} \int_{\mathcal{S}_\ell }f(x-u)\left({1\over 1-|u|^2}\right)^{1-{\alpha\over n}} du
\\\\ \ds
~\leq~\C\sum_{0\leq\ell<\sigma(x)} 2^{-\left(\sigma(x)-\ell\right)\left({1\over q}\right)} \Big(\M_\rho f\Big)^{p\over q}(x) \left\| f\right\|_{\L^p(\R^n)}^{1-{p\over q}}
\\\\ \ds
~\leq~\C_{p~q}~ \Big(\M_\rho f\Big)^{p\over q}(x) \left\| f\right\|_{\L^p(\R^n)}^{1-{p\over q}}.
\end{array}
\eeq
By putting together (\ref{Est1 Partial}) and (\ref{Est2 Partial}), we obtain (\ref{Regularity est}) as required.

 \section{Proof of Lemma Two}
 \setcounter{equation}{0}
We begin this section by recalling some well known estimates of Bessel functions. More discussions can be found in the book by Watson \cite{Watson}.
 
$\bullet$ For  $\mu>-{1\over 2}, \nu\in\R$ and $\rho>0$,  a Bessel function can be defined by   \bel{Bessel}
\J_{\mu+\i \nu}(\rho)~=~{(\rho/2)^{\mu+\i \nu}\over \pi^{1\over 2}\Gamma\left(\mu+{1\over 2}+\i \nu\right)} \int_{-1}^1 e^{\i \rho s} (1-s^2)^{\mu-{1\over 2}+\i \nu} ds.
\eeq
$\bullet$ For $\mu>-{1\over 2}, \nu\in\R$ and $\rho>0$, 
\bel{Bessel formula}
 \J_{\mu+\i \nu}(\rho)~=~\left({2\over \pi \rho}\right)^{1\over 2}\cos\left(\rho-{\pi\over 2}(\mu+\i \nu)-{\pi\over 4}\right)+\e_{\mu+\i \nu}(\rho)
\eeq
where
\bel{e asymptotic}
 |\e_{\mu+\i \nu}(\rho)|~\leq~\C_\mu~ e^{\hbox{\small{\bf c}}|\nu|}\left\{\begin{array}{lr}\ds \rho^{-{1\over 2}},\qquad 0<\rho\leq1,
 \\ \ds
 \rho^{-{3\over 2}},\qquad~~~ \rho>1.
 \end{array}\right.
 \eeq
 $\bullet$ For every $\mu, \nu\in\R$ and $\rho>0$, we have the identity
\bel{J identity}
\left({\mu+\i \nu\over \rho}\right)\J_{\mu+\i \nu}(\rho)~=~\J_{\mu-1+\i \nu}(\rho)+\J_{\mu+1+\i \nu}(\rho).
\eeq
All together, (\ref{Bessel}), (\ref{Bessel formula})-(\ref{e asymptotic}) and (\ref{J identity}) imply:

$\bullet$ For every $\mu, \nu\in\R$ and $\rho>0$,
\bel{J asymptotic}
\left| \rho^{-(\mu+\i \nu)}\J_{\mu+\i \nu}(\rho)\right|~\leq~\C_\mu~\left({1\over 1+\rho}\right)^{{1\over 2}+\mu}~e^{\c|\nu|}.
\eeq

 Let ${^s}\Omega^\alpha$ defined in (\ref{Omega s}) for $\Re\lambda(\alpha)+s<1$. We have
 \bel{Omega s norm compara}
 \left|{^s}\Omega^{n+\i\Im\alpha}(x)\right|~\leq~\C_s~e^{\hbox{\small{\bf c}}|\Im\alpha|}~\left({1\over 1-|x|^2}\right)^s,\qquad 0\leq s<1.
 \eeq
 By applying {\bf Lemma One} and then using (\ref{L^p omega}), we obtain
 \bel{Result Two < Est1}
 \begin{array}{cc}\ds
 \left\| f\ast\omega_s\ast{^s}\Omega^{n+\i\Im\alpha}\right\|_{\L^q(\R^n)}~\leq~\C_{p~q~s}~e^{\hbox{\small{\bf c}}|\Im\alpha|}~\left\| f\ast\omega_s\right\|_{\L^p(\R^n)}
 \\\\ \ds~~~~~~~~~~~~~~~~~~~~~~~~~~~~~~~~~~~~~~~~~~~~~~~~~~~~~~~~~~~~~~~~~~~
 ~\leq~\C_{p~q~s}~e^{\hbox{\small{\bf c}}|\Im\alpha|}~\left\| f\right\|_{\L^p(\R^n)},
 \qquad 1<p< q<\infty
 \\\\ \ds
 \hbox{ if}\qquad
 1-s~=~{1\over p}-{1\over q},\qquad 0<s<1.
 \end{array}
 \eeq
 From (\ref{Omega s norm compara}), it is also clear that
 \bel{Result Two < Est1 s=0}
  \left\| f\ast\Omega^{n+\i\Im\alpha}\right\|_{\L^\infty(\R^n)}~\leq~\C~e^{\hbox{\small{\bf c}}|\Im\alpha|}~\left\| f\right\|_{\L^1(\R^n)}.\qquad ( s=0 ) 
 \eeq
Let ${^s}\Omega^\alpha$ defined by its Fourier transform in (\ref{Omega s Transform}) for $s\ge0$ and $\alpha\in\Cx$. From (\ref{J asymptotic}), we have
 \bel{Omega s norm}
 \begin{array}{lr}\ds
 \left| {^s}\Hat{\Omega}^\alpha(\xi)\right|~\leq~\C_{\alpha~s}~e^{\hbox{\small{\bf c}}|\Im\alpha|}~\left({1\over 1+|\xi|}\right)^{\left({n+1\over 2n}\right)\Re\alpha-s}.
 \end{array}
 \eeq
 By using (\ref{Omega s norm}) and applying Plancherel theorem, we find
 \bel{L^2 Est}
 \begin{array}{lr}\ds
 \left\| f\ast\omega_s\ast{^s}\Omega^\alpha\right\|_{\L^2(\R^n)}~\leq~\C_s~e^{\hbox{\small{\bf c}}|\Im\alpha|}~\left\| f\ast\omega_s\right\|_{\L^2(\R^n)}
 \\\\ \ds ~~~~~~~~~~~~~~~~~~~~~~~~~~~~~~~~~
 ~\leq~\C_s~e^{\hbox{\small{\bf c}}|\Im\alpha|}~\left\| f\right\|_{\L^2(\R^n)}\qquad\hbox{\small{by (\ref{L^p omega}).}}
 \end{array}
 \eeq
 In particular, we have
  \bel{L^2 Est 0}
 \left\| f\ast\omega_s\ast{^s}\Omega^{0+\i\Im\alpha}\right\|_{\L^2(\R^n)} ~\leq~\C_s~e^{\hbox{\small{\bf c}}|\Im\alpha|}~\left\| f\right\|_{\L^2(\R^n)}.
 \eeq
Let $0<\Re\alpha<n$ and $0\leq s\leq {1\over 2}$. By using (\ref{Result Two < Est1})-(\ref{Result Two < Est1 s=0}) and (\ref{L^2 Est 0})  then applying Stein interpolation theorem \cite{Stein'}, we have
 \bel{Result Two < Est}
  \begin{array}{cc}\ds
 \left\| f\ast\omega_s\ast{^s}\Omega^{\alpha}\right\|_{\L^q(\R^n)}~\leq~\C_{p~q~s}~e^{\hbox{\small{\bf c}}|\Im\alpha|}~\left\| f\right\|_{\L^p(\R^n)}
 \\\\ \ds
 \hbox{ if}\qquad
 (1-s){\Re\alpha\over n}~=~{1\over p}-{1\over q},\qquad {1\over 2}+{\Re\alpha\over n}\left({1\over 2}-s\right)~\leq~{1\over p}~\leq~{1\over 2}+{\Re\alpha\over 2n}.
 \end{array}
 \eeq
  Let $0<\Re\alpha<n$ and $s\ge {1\over 2}$. From (\ref{Omega s norm}), we find
\bel{f Omega s}
\begin{array}{lr}\ds
\left|\Hat{f}(\xi)\Hat{\omega}_s(\xi){^s}\Hat{\Omega}^\alpha(\xi)\right|~\leq~\C_{\alpha~s}~e^{\hbox{\small{\bf c}}|\Im\alpha|}~\left|\Hat{f}(\xi)\right| \Hat{\omega}_s(\xi)\left(1+|\xi|\right)^{s} \left({1\over 1+|\xi|}\right)^{\left({n+1\over 2n}\right)\Re\alpha}
\\\\ \ds~~~~~~~~~~~~~~~~~~~~~~~~~~~~~
~\leq~\C_{\alpha~s}~e^{\hbox{\small{\bf c}}|\Im\alpha|}~\left|\Hat{f}(\xi)\right|\left({1\over |\xi|}\right)^{\Re\left({\alpha\over 2}\right)}.
\end{array}
\eeq
Note that the inverse Fourier transform of $|\xi|^{-\alpha}$ equals $\C_\alpha |x|^{\alpha-n}$ for $0<\alpha<n$. 

By applying Hardy-Littlewood-Sobolev inequality \cite{Hardy-Littlewood}-\cite{Sobolev} together with Plancherel theorem, we simultaneously have
  \bel{L^2 p q s result}
 \begin{array}{cc}\ds
  \left\| f\ast\omega_s\ast{^s}\Omega^\alpha\right\|_{\L^q(\R^n)}~\leq~\C_{q~s}~e^{\hbox{\small{\bf c}}|\Im\alpha|}~\left\|f\right\|_{\L^2(\R^n)}
\qquad
   \hbox{if}\qquad {\Re\alpha\over 2n}~=~{1\over 2}-{1\over q},
   \\\\ \ds
  \left\| f\ast\omega_s\ast{^s}\Omega^\alpha\right\|_{\L^2(\R^n)}~\leq~\C_{p~s}~e^{\hbox{\small{\bf c}}|\Im\alpha|}~\left\|f\right\|_{\L^p(\R^n)}
\qquad
   \hbox{if}\qquad {\Re\alpha\over 2n}~=~{1\over p}-{1\over 2}.
   \end{array}
   \eeq
 From (\ref{L^2 p q s result}), by applying Riesz-Thorin interpolation theorem, we find  
  \bel{Result Two > Est}
  \begin{array}{cc}\ds
 \left\| f\ast\omega_s\ast{^s}\Omega^{\alpha}\right\|_{\L^q(\R^n)}~\leq~\C_{p~q~s}~e^{\hbox{\small{\bf c}}|\Im\alpha|}~\left\| f\right\|_{\L^p(\R^n)}
 \\\\ \ds
 \hbox{ if}\qquad
{\Re\alpha\over 2n}~=~{1\over p}-{1\over q},\qquad {1\over 2}~\leq~{1\over p}~\leq~{1\over 2}+{\Re\alpha\over 2n}.
 \end{array}
 \eeq

 \section{Interpolation on a family of analytic operators}
 \setcounter{equation}{0}
Consider a family of analytic operators $f\ast\omega_s\ast\Theta^{\alpha~s}_z$ for $0<\Re z< 1$ where
\bel{Theta_z transform}
\begin{array}{lr}\ds
\Hat{\Theta}^{\alpha~s}_z(\xi)~\doteq~
\left({1\over|\xi|}\right)^{\left({n+1\over 2n}\right)\alpha-{1\over 2}+\left({n-1\over 2}\right)z-s(1-z)} 
\J_{\left({n+1\over 2n}\right)\alpha-{1\over 2}+\left({n-1\over 2}\right)z-s(1-z)}\Big(2\pi|\xi|\Big).
\end{array}
\eeq
In particular, we have
\bel{terms}
\begin{array}{cc}\ds
\Hat{\Theta}^{\alpha~s}_0(\xi)~=~{^\flat}\Hat{\Omega}^\alpha(\xi), \qquad \Hat{\Theta}^{\alpha~s}_1(\xi)~=~ {^s}\Hat{\Omega}^\alpha(\xi)
\\\\ \ds
\hbox{\small{and}}\qquad
 \Hat{\Theta}^{\alpha~s}_{2s\over n-1+2s}(\xi)~=~\left({1\over|\xi|}\right)^{\left({n+1\over 2n}\right)\alpha-{1\over 2}} \J_{\left({n+1\over 2n}\right)\alpha-{1\over 2}}\Big(2\pi|\xi|\Big)
  ~=~
 \Hat{\Omega}^\alpha(\xi)\qquad \hbox{\small{by (\ref{Omega^alpha Transform}).}}
 \end{array}
\eeq 
Let $0<\Re\alpha<{2n\over n+1}$. As a special case of (\ref{Omega flat regularity}), for $0\leq s\leq{1\over 2}$, we have
\bel{Omega flat regularity <}
\begin{array}{cc}\ds
\left\| f\ast\omega_s\ast{^\flat}\Omega^{\alpha}\right\|_{\L^q(\R^n)}~\leq~\C_{p~q}~e^{\hbox{\small{\bf c}}|\Im\alpha|}~\left\|f\ast\omega_s\right\|_{\L^p(\R^n)}
\\\\ \ds~~~~~~~~~~~~~~~~~~~~~~~~~~~~~~~~~~~~~~~~~~~~~~~
~\leq~\C_{p~q}~e^{\hbox{\small{\bf c}}|\Im\alpha|}~\left\|f\right\|_{\L^p(\R^n)}\qquad \hbox{\small{by (\ref{L^p omega})}}
\\\\ \ds~~~~
\hbox{if}\qquad \left({n+1\over 2n}\right)\Re\alpha~=~{1\over p}-{1\over q},\qquad 1<p<q<\infty.
\end{array}
\eeq
On the other hand, for $s\ge{1\over 2}$, we have
\bel{Omega flat regularity >}
\begin{array}{cc}\ds
\left\| f\ast\omega_s\ast{^\flat}\Omega^{\alpha}\right\|_{\L^q(\R^n)}~\leq~\C_{p~q}~e^{\hbox{\small{\bf c}}|\Im\alpha|}~\left\|f\ast\omega_s\right\|_{\L^p(\R^n)}
\\\\ \ds~~~~~~~~~~~~~~~~~~~~~~~~~~~~~~~~~~~~~~~~~~~~~~~~~
~\leq~\C_{p~q}~e^{\hbox{\small{\bf c}}|\Im\alpha|}~\left\|f\right\|_{\L^p(\R^n)},\qquad \hbox{\small{by (\ref{L^p omega})}}
\\\\ \ds
\hbox{if}\qquad \left({n-1+4s\over 4s}\right){\Re\alpha\over n}~=~{1\over p}-{1\over q},\qquad 1<p<q<\infty.~~~~
\end{array}
\eeq
Recall {\bf Lemma Two}.  Let $0< s\leq{1\over 2}$. By using (\ref{Result Two < EST}) and (\ref{Omega flat regularity <}) respectively for $f\ast\omega_s\ast\Theta^{\alpha~s}_z$  defined as (\ref{Theta_z transform})-(\ref{terms}) and 
 applying Stein interpolation theorem \cite{Stein'}, we have
\bel{EST <}
\begin{array}{lr}\ds
~~~~~~~~~~~~~ \left\| f\ast\omega_s\ast\Omega^{\alpha}\right\|_{\L^q(\R^n)}~\leq~\C_{p~q}~e^{\hbox{\small{\bf c}}|\Im\alpha|}~\left\| f\right\|_{\L^p(\R^n)},\qquad 1<p<q<\infty
 \\\\ \ds
 \hbox{if}\qquad {1\over p}-{1\over q}~=~\left({2s\over n-1+2s}\right)\left({n+1\over 2n}\right)\Re\alpha+\left({n-1\over n-1+2s}\right)\left({1-s\over n}\right)\Re\alpha
 ~=~{\Re\alpha\over n}
 \\\\ \ds
 \hbox{and}\qquad \left({2s\over n-1+2s}\right) \left({n+1\over 2n}\right)\Re\alpha+ \left({n-1\over n-1+2s}\right)\left[ {1\over 2}+{\Re\alpha\over n}\left({1\over 2}-s\right)\right] 
 \\\\ \ds~~~~~~~~~~~~~~~~
 ~=~{n-1\over 2n-2+4s}+\left({4s+n-1\over 2n-2+4s}\right){\Re\alpha\over n}~<~{1\over p}~<~
 \\\\ \ds
  \left({2s\over n-1+2s}\right)+ \left({n-1\over n-1+2s}\right) \left[{1\over 2}+{\Re\alpha\over 2n}\right] 
 ~=~{n-1+4s\over 2n-2+4s}+\left({n-1\over 2n-2+4s}\right){\Re\alpha\over n}. 
 \end{array}
 \eeq
 
Let $s\ge{1\over 2}$. By using (\ref{Result Two > EST}) and (\ref{Omega flat regularity >}) respectively for $f\ast\omega_s\ast\Theta^{\alpha~s}_z$  defined as (\ref{Theta_z transform})-(\ref{terms}) and 
 applying Stein interpolation theorem \cite{Stein'}, we have
 \bel{EST >}
\begin{array}{lr}\ds
~~~~~~~~~~~~~ \left\| f\ast\omega_s\ast\Omega^{\alpha}\right\|_{\L^q(\R^n)}~\leq~\C_{p~q}~e^{\hbox{\small{\bf c}}|\Im\alpha|}~\left\| f\right\|_{\L^p(\R^n)},\qquad 1<p<q<\infty
 \\\\ \ds
 \hbox{if}\qquad {1\over p}-{1\over q}~=~\left({2s\over n-1+2s}\right)\left({n-1+4s\over 4s}\right){\Re\alpha\over n}+\left({n-1\over n-1+2s}\right){\Re\alpha\over 2n}
 ~=~{\Re\alpha\over n}
 \\\\ \ds
 \hbox{and}\qquad \left({2s\over n-1+2s}\right) \left({n-1+4s\over 4s}\right){\Re\alpha\over n}+ \left({n-1\over n-1+2s}\right){1\over 2}
 \\\\ \ds~~~~~~~~~~~~~~~~
 ~=~{n-1\over 2n-2+4s}+\left({4s+n-1\over 2n-2+4s}\right){\Re\alpha\over n}~<~{1\over p}~<~
 \\\\ \ds
  \left({2s\over n-1+2s}\right)+ \left({n-1\over n-1+2s}\right) \left[{1\over 2}+{\Re\alpha\over 2n}\right] 
 ~=~{n-1+4s\over 2n-2+4s}+\left({n-1\over 2n-2+4s}\right){\Re\alpha\over n}. 
 \end{array}
 \eeq
By putting together (\ref{EST <}) and (\ref{EST >}), we conclude
\bel{EST}
\begin{array}{ccc}\ds
 \left\| f\ast\omega_s\ast\Omega^{\alpha}\right\|_{\L^q(\R^n)}~\leq~\C_{p~q}~e^{\hbox{\small{\bf c}}|\Im\alpha|}~\left\| f\right\|_{\L^p(\R^n)}
\\\\ \ds
\hbox{if}\qquad {\Re\alpha\over n}~=~{1\over p}-{1\over q}\qquad \hbox{and}
\\\\ \ds
{n-1\over 2n-2+4s}+\left({4s+n-1\over 2n-2+4s}\right){\Re\alpha\over n}~<~{1\over p}~<~{n-1+4s\over 2n-2+4s}+\left({n-1\over 2n-2+4s}\right){\Re\alpha\over n}
\end{array}
\eeq
for $s>0$ and $0<\Re\alpha<{2n\over n+1}$.

On the other hand, from (\ref{Result Two < Est}) and (\ref{L^p omega}), we have
 \bel{Result Two < Est s=0}
  \begin{array}{cc}\ds
 \left\| f\ast\omega_s\ast\Omega^{\alpha}\right\|_{\L^q(\R^n)}~\leq~\C_{p~q}~e^{\hbox{\small{\bf c}}|\Im\alpha|}~\left\| f\ast\omega_s\right\|_{\L^p(\R^n)}
 \\\\ \ds~~~~~~~~~~~~~~~~~~~~~~~~
 ~\leq~\C_{p~q}~e^{\hbox{\small{\bf c}}|\Im\alpha|}~\left\| f\right\|_{\L^p(\R^n)}
 \\\\ \ds
 \hbox{ for}\qquad
 {1\over p}~=~{1\over 2}+{\Re\alpha\over 2n},\qquad {1\over q}~=~{1\over 2}-{\Re\alpha\over 2n} \end{array}
 \eeq
 and $0<\Re\alpha<n$.

By using (\ref{EST})-(\ref{Result Two < Est s=0}) and applying Stein interpolation theorem \cite{Stein'}, we find
\bel{RESULT ONE s}
\begin{array}{ccc}\ds
 \left\| f\ast\omega_s\ast\Omega^{\alpha}\right\|_{\L^q(\R^n)}~\leq~\C_{p~q}~\left\| f\right\|_{\L^p(\R^n)},\qquad s>0
\\\\ \ds
\hbox{if}\qquad {\alpha\over n}~=~{1\over p}-{1\over q}\qquad \hbox{and}
\\\\ \ds
{n-1\over 2n-2+4s}+\left({4s+n-1\over 2n-2+4s}\right){\alpha\over n}~<~{1\over p}~<~{n-1+4s\over 2n-2+4s}+\left({n-1\over 2n-2+4s}\right){\alpha\over n}.
\end{array}
\eeq
Note that $f\ast\Omega^\alpha=\F\ast\omega_s\ast\Omega^\alpha$ of which $\Hat{\F}(\xi)\doteq\Hat{f}(\xi)\left(1+|\xi|^2\right)^{s\over 2}$. We have $\F\in\L^p(\R^n)$ provided that $f\in\L^p_s(\R^n)$. This completes the proof of {\bf Theorem One}.

Let $0<\Re\alpha<{2n\over n+1}$. As a special case of (\ref{Omega flat regularity}), for every $ s\ge0$, we have
\bel{Omega flat regularity p=q}
\begin{array}{lr}\ds
\left\| f\ast\omega_s\ast{^\flat}\Omega^{\alpha}\right\|_{\L^p(\R^n)}~\leq~\C_{p~\alpha}~e^{\hbox{\small{\bf c}}|\Im\alpha|}~\left\|f\ast\omega_s\right\|_{\L^p(\R^n)}
\\\\ \ds~~~~~~~~~~~~~~~~~~~~~~~~~~~~~~~~~
~\leq~\C_{p~\alpha}~e^{\hbox{\small{\bf c}}|\Im\alpha|}~\left\|f\right\|_{\L^p(\R^n)},\qquad 1<p<\infty\qquad \hbox{\small{by (\ref{L^p omega}).}}
\end{array}
\eeq
Recall {\bf Lemma Two}. Let $s>0$. By using (\ref{Result Two L^2}) and (\ref{Omega flat regularity p=q}) respectively for $f\ast\omega_s\ast\Theta^{\alpha~s}_z$  defined as (\ref{Theta_z transform})-(\ref{terms}) and 
 applying Stein interpolation theorem \cite{Stein'}, we have
\bel{EST small}
\begin{array}{cc}\ds
\left\| f\ast\omega_s\ast\Omega^{\alpha}\right\|_{\L^p(\R^n)}~\leq~\C_{p~\alpha}~e^{\hbox{\small{\bf c}}|\Im\alpha|}~\left\| f\right\|_{\L^p(\R^n)}
 \\\\ \ds
 \hbox{for}\qquad \left({n-1\over n-1+2s}\right){1\over 2}+0 ~=~
 {n-1\over 2n-2+4s}
 ~<~{1\over p}~<~
 \\\\ \ds~~~~~~~
 \left({n-1\over n-1+2s}\right){1\over 2}+{2s\over n-1+2s} ~=~{n-1+4s\over 2n-2+4s}.
 \end{array}
 \eeq
Let $\Omega^\alpha$ defined in (\ref{Omega^alpha}) for $\Re\lambda(\alpha)={n+1\over 2}\left(1-{\Re\alpha\over n}\right)<1$. Observe that
\bel{Omega Compara Size}
\left|\Omega^\alpha(x)\right|~\leq~\pi^{-\Re\lambda(\alpha)} \left|\Gamma^{-1}(1-\lambda(\alpha))\right|{^\natural}\Omega^{\left({n+1\over 2}\right)\Re\alpha-\left({n-1\over 2}\right)n}(x)
\eeq
for $\left({n-1\over n+1}\right)n<\Re\alpha<n$ where ${^\natural}\Omega^\alpha$ is defined in (\ref{Omega^alpha natural}) for $0<\alpha<n$. 

By using (\ref{Omega Compara Size}) and applying {\bf Lemma One}, we have
\bel{L^p EST large}
\begin{array}{lr}\ds
\left\| f\ast\omega_s\ast\Omega^\alpha\right\|_{\L^p(\R^n)} ~\leq~\C_{p~\alpha}~e^{\hbox{\small{\bf c}}|\Im\alpha|}~\left\| f\ast\omega_s\right\|_{\L^p(\R^n)}
\\\\ \ds~~~~~~~~~~~~~~~~~~~~~~~~~~~~~~~
~\leq~\C_{p~\alpha}~e^{\hbox{\small{\bf c}}|\Im\alpha|}~\left\| f\right\|_{\L^p(\R^n)},\qquad 1<p<\infty\qquad\hbox{\small{by (\ref{L^p omega})}}
\end{array}
\eeq
whenever $\left({n-1\over n+1}\right)n<\Re\alpha<n$. Note that $\Hat{\Omega}^\alpha(\xi)$ at $\alpha=\left({n-1\over n+1}\right)n$ is the Fourier transform of the surface measure on $\mathds{S}^{n-1}$.  We therefore conclude {\bf Remark One}.

By using (\ref{EST small}) with (\ref{L^p EST large}) and applying Stein interpolation theorem \cite{Stein'}, we have
\bel{EST L^p}
\begin{array}{ccc}\ds
\left\| f\ast\omega_s\ast\Omega^\alpha\right\|_{\L^p(\R^n)}
 ~\leq~\C_{p~\alpha}~\left\| f\right\|_{\L^p(\R^n)}
 \\\\ \ds
 \hbox{for}\qquad \left[1-{\alpha\over n}\left({n+1\over n-1}\right)\right] {n-1\over 2n-2+4s}~=~{n-1\over 2n-2+4s}-{\alpha\over n}\left({n+1\over 2n-2+4s}\right)
 ~<~{1\over p}~<~
 \\\\ \ds~~~~~~~~~~~~~~~
 \left[1-{\alpha\over n}\left({n+1\over n-1}\right)\right] {n-1+4s\over 2n-2+4s} +{\alpha\over n}\left({n+1\over n-1}\right)~=~{n-1+4s\over 2n-2+4s}+{\alpha\over n}\left({n+1\over 2n-2+4s}\right)
 \end{array}
\eeq
and $0<\alpha< \left({n-1\over n+1}\right)n$.

Lastly, $f\ast\Omega^\alpha=\F\ast\omega_s\ast\Omega^\alpha$ where $\Hat{\F}(\xi)=\Hat{f}(\xi)\left(1+|\xi|^2\right)^{s\over 2}$. Moreover, $\F\in\L^p(\R^n)$ because $f\in\L^p_s(\R^n)$. We finish the proof of  {\bf Theorem Two}.

\v

\small{Department of Mathematics, Westlake University.} 

 \small{email: wangzipeng@westlake.edu.cn}


\begin{thebibliography}{100}
\bibitem{Hardy-Littlewood}
{\small G.~H.~Hardy and J.~E.~Littlewood, {\it Some Properties of Fractional Integrals}, 
Mathematische Zeitschrift {\bf 27}: 565-606, 1928}.







\bibitem{Sobolev}
\small{ S.~L.~Sobolev, {\it On a Theorem of Functional Analysis}, Matematicheskii Sbornik {\bf 46}: 471-497, 1938}.



\bibitem{Strichartz}{\small R. Strichartz, {\it Convolutions with kernels having singularity on a sphere}, Transaction of the American Mathematical Society {\bf 148}: 461-471, 1970}.

\bibitem{Strichartz'}R. Strichartz, {\it A priori estimates for the wave equation and some applications}, Journal of Functional Analysis {\bf 5}: 218-235, 1970.



\bibitem{Hedberg}{\small L.~Hedberg, {\it On Certain Convolution Inequalities}, Proceeding of American Mathematical Society {\bf 36}: 505-510, 1972}.





  \bibitem{Stein'}{\small E.~M.~Stein, {\it Interpolation of linear operators}, Transaction of the American Mathematical Society {\bf 83}: 482-492, 1956.}



\bibitem{Stein}{\small E.~M.~Stein, {\it
Harmonic Analysis: Real-Variable Methods, Orthogonality and Oscillatory Integrals},
 Princeton University Press, 1993}.
 
 

 
 
 
 
 
\bibitem{Stein-Weiss*}{\small E.~M.~Stein and G.~Weiss, {\it Introduction to Fourier Analysis on Euclidean Space},
 Princeton University Press, 1971}.
 
 
 
 
 

\bibitem{Fefferman}{\small C.~Fefferman, {\it A note on Spherical Summation Multipliers}, Israel Journal of Mathematics, {\bf 15}: 44-52, 1973.}



\bibitem{Sogge}{\small C.~D.~Sogge, {\it $\L^p$ estimates for the wave equation and applications}, Journ\'{e}es \'{E}quations aux d\'{e}riv\'{e}es partielles, 1-12, 1993.}


\bibitem{Sogge'}{\small C.~D.~Sogge, {\it Lectures on Non-linear Wave Equations} (second edition), Johns Hopkins University Press, 2013.}





\bibitem{Watson}{\small G.~N.~Watson, {\it Theory of Bessel Functions}, Cambridge University Press, Cambridge, 1944}.




\end{thebibliography}
\end{document}